\documentclass[a4paper,12pt]{amsart}
\usepackage{pst-all, pstricks, pst-node}
\usepackage{amssymb, amsmath, amsthm, amscd}
\usepackage{subfigure}
\usepackage{graphicx, hhline, a4wide, ulem}
\usepackage[arrow,tips,matrix]{xy}
\usepackage[arrow,matrix]{xy}
\usepackage{auto-pst-pdf}

\theoremstyle{ams}
\newtheorem{theorem}{Theorem}[section]
\newtheorem{proposition}[theorem]{Proposition}
\newtheorem{lemma}[theorem]{Lemma}
\newtheorem{corollary}[theorem]{Corollary}

\theoremstyle{definition}

\newtheorem{definition}[theorem]{Definition}
\newtheorem{remark}[theorem]{Remark}

\numberwithin{table}{section} \numberwithin{figure}{section}
\numberwithin{equation}{section}

\newcommand{\C}{\mathbb{C}}

\newcommand{\N}{\mathbb{N}}
\newcommand{\Q}{\mathbb{Q}}
\newcommand{\R}{\mathbb{R}}
\newcommand{\Z}{\mathbb{Z}}

\begin{document}
\title[Embedded surfaces for symplectic circle actions]{Embedded surfaces for symplectic circle actions}

\author[Y. Cho]{Yunhyung Cho}
\address{Center for Geometry and Physics, Institute for Basic Science(IBS), Pohang, Republic of Korea 37673}
\email{yhcho@ibs.re.kr}
\author[M. K. Kim]{Min Kyu Kim}
\address{Department of Mathematics Education,
Gyeongin National University of Education, San 59-12, Gyesan-dong,
Gyeyang-gu, Incheon, 407-753, Korea}
\email{mkkim@kias.re.kr}
\author[D. Y. Suh]{Dong Youp Suh}
\address{Department of Mathematical Sciences, KAIST, 335 Gwahangno, Yu-sung Gu, Daejeon 305-701, Korea}
\email{dysuh@math.kaist.ac.kr}

\thanks{The first author was supported by IBS-R003-D1. The second author is partially supported by Gyeongin National University of Education research fund. The third author was supported in part by Basic Science Research Program through the National Research Foundation of Korea(NRF) grant funded by the Ministry of Education(2013R1A1A2007780)}

\date{\today}
\maketitle

\begin{abstract}

     The purpose of this article is to characterize symplectic and Hamiltonian circle actions on symplectic manifolds in terms of symplectic embeddings of Riemann surfaces.
     More precisely, we will show that (1) if $(M,\omega)$ admits a Hamiltonian $S^1$-action, then there exists an $S^1$-invariant symplectic $2$-sphere $S$ in $(M,\omega)$ such that $\langle c_1(M), [S] \rangle > 0$, and (2) if the action is non-Hamiltonian, then there exists an $S^1$-invariant symplectic $2$-torus $T$ in $(M,\omega)$ such that $\langle c_1(M), [T] \rangle = 0$.
     As applications, we will give a very simple proof of the following well-known theorem which was proved by Atiyah-Bott \cite{AB}, Lupton-Oprea \cite{LO}, and Ono \cite{O2} : suppose that $(M,\omega)$ is a smooth closed symplectic manifold satisfying $c_1(TM)=\lambda \cdot [\omega]$ for some $\lambda \in \R$ and let $G$ be a compact connected Lie group acting effectively on $M$ preserving $\omega$. Then (1) if $\lambda < 0$, then $G$ must be trivial, (2) if $\lambda=0$, then the $G$-action is non-Hamiltonian, and (3) if $\lambda > 0$, then the $G$-action is Hamiltonian.

\end{abstract}

\section{Introduction}

%This expository article is to give a complete proof of the following theorem.
%\begin{theorem}\label{theorem : main}
%    Let $(M,\omega)$ be a smooth closed symplectic manifold such that $c_1(TM,J) = \lambda[\omega] \in H^2(M;\Z)$ for an $\omega$-compatible almost complex structure $J$ on $M$.
%    If there is a compact connected Lie group $G$ acting effectively on $M$ preserving $\omega$, then
%    \begin{enumerate}
%        \item if $\lambda > 0$, then the $G$-action must be Hamiltonian,
%        \item if $\lambda = 0$, then the $G$-action must be non-Hamiltonian, and
%        \item if $\lambda < 0$, then $G$ must be trivial.
%    \end{enumerate}
%\end{theorem}
The purpose of this article is to characterize symplectic and Hamiltonian circle actions on symplectic manifolds in terms of symplectic embeddings of Riemann surfaces.
We first consider the following simple situation which provides the motivation of our work.

Let $\Sigma_g$ be a two-dimensional smooth closed oriented manifold with genus $g$, and let $\mathrm{Diff}(\Sigma_g)$ denote the \textit{diffeomorphism group}
of $\Sigma_g$.
Suppose $G$ is a compact connected Lie group acting on $\Sigma_g$ effectively\footnote{A
$G$-action on a manifold $M$ is called \textit{effective} if the identity element $1 \in G$ is the unique element
which fixes whole $M$.}, i.e., there is an injective Lie group homomorphism $$ \phi : G \hookrightarrow
\mathrm{Diff}(\Sigma_g).$$ Since $G$ is compact and connected, by averaging any given Riemannian metric over the Haar measure of
$G$, we can get a $G$-invariant metric $\Omega$ on $\Sigma_g$ so that we may regard $G$ as a closed Lie
subgroup of the identity component $\mathrm{Iso}(\Sigma_g, \Omega)_0$ of the full isometry group
$\mathrm{Iso}(\Sigma_g, \Omega)$. In particular, we have the following inequalities :
$$ \dim G \leq \dim \mathrm{Iso}(\Sigma_g,\Omega)_0 = \dim \mathrm{Iso}(\Sigma_g,\Omega) \leq 3 $$ where the last inequality comes from the classical fact
that $$\dim \mathrm{Iso}(M^n, h) \leq \frac{n(n+1)}{2}$$ for any $n$-dimensional complete Riemannian manifold
$(M^n,h)$, see \cite{Ma} for more details. Note that the isometry groups of a 2-sphere $S^2$ and a 2-torus $T^2$ are
well-known such that
\begin{itemize}
    \item $\mathrm{Iso}(S^2,h)_0$ is a subgroup of $\mathrm{SO}(3)$ for any metric $h$ on $S^2$, and
    \item $\mathrm{Iso}(T^2,h)_0$ is a subgroup of $\mathrm{SO}(2)\times \mathrm{SO}(2)$ for any metric $h$ on $T^2$.
\end{itemize}
In the case when the genus $g \geq 2$, then the isometry
group of $\Sigma_g$ is finite with respect to any metric on
$\Sigma_g$, see \cite[Chapter 7]{FM} for more details.Thus the condition $g \geq 2$ is an obstruction to the existence of an action of a compact Lie group of positive dimension on Riemann surfaces.
Together with the fact that $\chi(\Sigma_g) = 2-2g$, one might expect that the existence of a non-trivial compact connected Lie group $G$ in $\mathrm{Diff}(M)$ might be obstructed by the Euler characteristic $\chi(M)$, i.e., the negativity of $\chi(M)$ would imply the non-existence of a $G$-action such as the case of compact Riemann surfaces.
Unfortunately, it is inappropriate to use $\chi(M)$ since if we let $M = S^2 \times \Sigma_g$ for $g \geq 2$, then $M$ admits a circle action on the first factor but $M$ satisfies $\chi(M) < 0$.
Therefore, instead of $\chi$, we use some geometric structure as follows.
Let us consider a $G$-invariant volume form $\omega$ on $\Sigma_g$.
Then $\omega$ is a symplectic form\footnote{A differential two-form $\omega$ on a manifold $M$ is called \textit{symplectic} if it is closed ($d\omega = 0$) and non-degenerate ($\omega_p : T_pM \times T_pM \rightarrow \R$ is a non-degenerate bilinear form for every $p \in M$).} and thus $(\Sigma_g, \omega)$ is a \textit{symplectic manifold}\footnote{A symplectic manifold is a pair $(M,\omega)$ which consists of a smooth manifold $M$ and a symplectic form $\omega$ on $M$}.
Hence we may think of $G$ as a subgroup of the \textit{symplectomorphism group}
$$ ~\mathrm{Symp}(\Sigma_g, \omega) := \{ g\in \mathrm{Diff}(\Sigma_g) ~|~ g^*\omega = \omega \}. $$
Let $J$ be an $\omega$-compatible almost complex structure
on $\Sigma_g$, i.e., $\omega(\cdot, \cdot) = \omega(J\cdot, J\cdot)$ and $\omega(\cdot,
J\cdot)$ is a Riemannian metric on $M$. Note that for any symplectic manifold $(M,\omega)$, $\omega$-compatible almost complex structure $J$ always exists.
In fact, the space $\mathcal{J}(M,\omega)$ of $\omega$-compatible almost complex structures is a contractible space (see \cite{McS}), which implies that $(TM,J)$ and $(TM,J')$ are
isomorphic as a complex vector bundle for any $J$ and $J'$ in $\mathcal{J}(M,\omega)$.
Hence the first Chern class $c_1(M,\omega)$ of $(TM,J)$ does not depend on the choice of a $\omega$-compatible almost complex structure $J$.
We denote by $c_1(M,\omega)$ the first Chern class of $(TM,J)$ for $J \in \mathcal{J}(M,\omega)$.

Since $\omega$ represents a non-zero element $[\omega] \in H^2(\Sigma_g;\R) \cong \R$, there exists a constant $\lambda \in \R$ such that
$$ c_1(\Sigma_g, \omega) = \lambda \cdot [\omega].$$
Also, the Chern number $\langle c_1(\Sigma_g,\omega), [\Sigma_g] \rangle$ is the same as the Euler characteristic $\chi(\Sigma_g) = 2 - 2g$ of $\Sigma_g$ where $[\Sigma_g] \in
H_2(\Sigma_g,\Z)$ is the fundamental homology class of $\Sigma_g$. Hence we may conclude as follows.
\begin{itemize}
    \item If $(\Sigma_g, \omega)$ is a closed symplectic surface such that $\langle c_1(\Sigma_g,\omega), [\Sigma_g] \rangle < 0$, then there is no compact connected Lie group action preserving $\omega$.
\end{itemize}
In fact, we can say more about the $G$-action on $(\Sigma_g, \omega)$ in the symplectic setting.
For a given symplectic manifold $(M,\omega)$ with the symplectomorphism group $\mathrm{Symp}(M,\omega)$, we say that $\phi \in \mathrm{Symp}(M,\omega)$ is a \textit{Hamiltonian diffeomorphism} if there exists an isotopy $\phi_t$ for $0 \leq t \leq 1$ such that
\begin{itemize}
    \item $\phi_0 = \mathrm{id}$,
    \item $\phi_1 = \phi$, and
    \item $i_{X_t}\omega = \omega(X_t, \cdot) = dH_t$ for some family of smooth functions $\{H_t : M \rightarrow \R \}_{0 \leq t \leq 1}$ where $X_t$ is a time-dependent vector field on $M$ such that
    \begin{equation}\label{equation : flow equation}
         X_t \circ \phi_t = \frac{d}{dt}\phi_t.
    \end{equation}
\end{itemize}
We denote by $\mathrm{Ham}(M,\omega)$ the set of all Hamiltonian diffeomorphisms on $(M,\omega)$ and it forms a group under the composition. In fact, the Hamiltonian diffeomorphism group
$\mathrm{Ham}(M,\omega)$ is path-connected and it is a normal subgroup of $\mathrm{Symp}(M,\omega)$, see \cite{McS} for more details.

Now suppose that $G$ is a compact connected Lie group acting on $(M,\omega)$ effectively and $\mathfrak{g}$ is the Lie algebra of $G$.
By definition of $\mathrm{Symp}(M,\omega)$ and by the surjectivity of the exponential map $\mathrm{exp} : T_eG \rightarrow G$,
we can easily show that $G$ is a subgroup of
$\mathrm{Symp}(M,\omega)$ if and only if each one-parameter subgroup generated by each element $X \in \mathfrak{g}$ preserves $\omega$, i.e.,
$$\mathcal{L}_{\underline{X}}\omega = (i_{\underline{X}}
\circ d + d \circ i_{\underline{X}})\omega = d \circ
i_{\underline{X}}\omega = 0$$ for every $X \in T_eG$ where
$\underline{X}$ is the vector field generated by $X$, i.e., $$\underline{X}_p := \frac{d}{dt}|_{t=0} (\exp tX) \cdot p$$ for each $p \in M$. In particular, if $i_{\underline{X}}\omega$ is exact for every $X \in \mathfrak{g}$, then $i_{\underline{X}}\omega = dH_X$ for some smooth function $H_X : M \rightarrow \R$ for each $X \in \mathfrak{g}$. Then, the family $\{\exp tX \}_{0 \leq t \leq 1}$ is an isotopy which connects $\mathrm{id} = \exp (0\cdot X)$ with $\exp (1\cdot X)$, and it satisfies the equation (\ref{equation : flow equation}), which means that $\exp X \in \mathrm{Ham}(M,\omega)$. By surjectivity of $\exp : \mathfrak{g} \rightarrow G$ (since $G$ is compact), we can deduce that $G$ is a subgroup of
$\mathrm{Ham}(M,\omega)$ if $i_{\underline{X}}\omega$ is exact for every $X \in T_eG$. Conversely, if $G$ is a subgroup of $\mathrm{Ham}(M,\omega)$, then one can easily check that $i_{\underline{X}}\omega$ is exact for every $X \in \mathfrak{g}$.
\begin{definition}\label{definition : Hamiltonian action}
    Let $G$ be a compact connected Lie group acting on $(M,\omega)$. We say that a $G$-action is \textit{Hamiltonian} if $i_{\underline{X}}\omega$ is exact for every $X \in \mathfrak{g}$.
    Equivalently, a $G$-action is Hamiltonian if $G$ acts on $M$ as a subgroup of $\mathrm{Ham}(M,\omega)$.
\end{definition}
Hence if the $G$-action is Hamiltonian, there exists a smooth map $H : M \rightarrow \mathfrak{g}^*$ (called a \textit{moment map}) such that $\langle H, X \rangle = H_X$ with $i_{\underline{X}}\omega = dH_X$ for every $X \in \mathfrak{g}$, where $\langle, \rangle : \mathfrak{g}^* \times \mathfrak{g} \rightarrow \R$ is the usual pairing of the Lie algebra $\mathfrak{g}$ with its dual.
Note that for each $X \in \mathfrak{g}$, the set of all critical points of $H_X$ coincides with the set of all points fixed by the
one-parameter subgroup generated by $X$ by the non-degeneracy of $\omega$. Note
that if $M$ is compact, then any smooth function on $M$ has
at least two critical points attaining its extrema. Hence if the $G$-action is Hamiltonian,
the one-parameter subgroup $\{\mathrm{exp}(tX) \}_{t \in
\R}$ has at least two fixed points for every $X \in \mathfrak{g}$. Therefore, we
have the following proposition as follows.

\begin{proposition}\label{proposition : two dimensional theorem}
    Let $(\Sigma_g,\omega)$ be a closed two-dimensional symplectic manifold with genus $g$ such that $c_1(\Sigma_g,\omega) = \lambda \cdot [\omega]$ and let $G$ be a compact connected Lie group. Suppose that the $G$ action is effective and it preserves $\omega$. Then,
\begin{enumerate}
    \item if $\lambda > 0$ ($g=0$), then the $G$-action is Hamiltonian.
    \item if $\lambda = 0$ ($g=1$), then the $G$-action is non-Hamiltonian, and
    \item if $\lambda < 0$ ($g > 1$), then $G = \{ 1 \}$.
\end{enumerate}
\end{proposition}

\begin{proof}
        If $g=0$, then $\Sigma_0 \cong S^2$ is simply connected. In particular, we have $H^1(\Sigma_0 ; \R) = 0$ so that $di_{\underline{X}}\omega = 0$ if and only if $i_{\underline{X}}\omega$ is exact, which implies that any symplectic $G$-action on $\Sigma_0$ is automatically Hamiltonian. For the second statement, recall that $SO(2) \times SO(2)$ acts on $\Sigma_1$ freely and hence the $G$-action on $\Sigma_1$ is also free. In particular, every one-parameter subgroup of $G$ has no fixed point so that the action is non-Hamiltonian. The last statement comes from the fact that the isometry group of $\Sigma_g$ with $g \geq 2$ is finite (see Page 2).
\end{proof}

In this point of view, we may think of the generalization of the closed Riemann surface case to the symplectic category.
Note that if $(M,\omega)$ is a symplectic manifold such that $c_1(M,\omega) = \lambda \cdot [\omega]$, then we call $(M,\omega)$ a \textit{monotone symplectic manifold}
if $\lambda > 0$, a \textit{symplectic Calabi-Yau manifold} if $\lambda = 0$, and a \textit{negatively monotone symplectic manifold} if
$\lambda < 0$. In this article, we will regard those three
families of symplectic manifolds as a generalization of
closed Riemann surfaces and we prove the following theorem, which is a generalization of Proposition \ref{proposition : two dimensional theorem} to the higher dimensional cases.

\begin{theorem}\label{theorem : main}
    Let $(M,\omega)$ be any smooth closed symplectic manifold such that $c_1(M,\omega) = \lambda \cdot [\omega]$ for some $\lambda \in \R$. Let $G$ be a compact connected Lie group which acts on $(M,\omega)$ effectively and preserves $\omega$. Then,
    \begin{enumerate}
        \item If $\lambda > 0$, then the $G$-action is Hamiltonian.
        \item If $\lambda = 0$, then the $G$-action is non-Hamiltonian.
        \item If $\lambda < 0$, then $G$ is trivial.
    \end{enumerate}
\end{theorem}

Note that Theorem \ref{theorem : main} is not new. Theorem \ref{theorem : main} (1) was already proved independently by Atiyah-Bott \cite{AB} and Lupton-Oprea \cite{LO}.
Also Theorem \ref{theorem : main} (2), (3) was proved by Ono \cite{O}. The proofs given by Atiyah-Bott \cite{AB} and Ono \cite{O} are based on the equivariant cohomology theory, and the proof of Lupton-Oprea \cite{LO} is based on the homotopy theory, in particular the theory of Gottlieb groups. We do not refer to the details of their proofs, instead we give much more elementary and simple proof of Theorem \ref{theorem : main} which can be obtained as a corollary of the following series of the propositions.

\begin{proposition}\label{proposition : two torus embedding}
    Let $(M,\omega)$ be a smooth closed symplectic manifold equipped with a smooth $S^1$-action preserving $\omega$.
    Suppose that $[\omega]$ is a rational class in $H^2(M;\Q)$. Then the action is non-Hamiltonian if and only if there exists an $S^1$-equivariant symplectic
    embedding of 2-torus $i \colon T^2 \hookrightarrow M$, where the circle acts freely on the left factor of $T^2 \cong S^1 \times S^1$.
    In particular, the normal bundle of $i(T^2)$ in $M$ is trivial so that $\langle c_1(M,\omega), [i(T^2)] \rangle = 0$.
\end{proposition}

\begin{proposition}\label{proposition : embedding sphere Hamiltonian case}
    Suppose that $(M,\omega)$ is a smooth closed symplectic manifold equipped with a Hamiltonian circle action.
    Then there exist an $S^1$-equivariant symplectic embedding of two-sphere $S$ satisfying $\langle c_1(TM), [S] \rangle > 0$.
\end{proposition}

\begin{proof}[Proof of Theorem \ref{theorem : main}]
First, suppose that the action is Hamiltonian. If $G$ is not trivial, then there exists a maximal subtorus $T$ with positive dimension in $G$. For any choice of a circle subgroup $S^1$ of $T$, there exists a symplectic two sphere $S$ in $M$ such that $ \langle c_1(TM), [S] \rangle = \lambda \cdot \langle [\omega], [S] \rangle > 0$ by Proposition \ref{proposition : embedding sphere Hamiltonian case}.
Since $\langle [\omega], [S] \rangle$ is a symplectic area of $S$, $\lambda$ must be positive.
Secondly, if the action is non-Hamiltonian, we can easily show that it is non-Hamiltonian with respect to $k \cdot  \omega$ for any nonzero $k \in \R$. 
If $c_1(M) \neq 0$, i.e., $\lambda \neq 0$, then the action is non-Hamiltonian with respect to $\lambda \cdot \omega$ which represents the integral class $c_1(M,\omega)$.
By Proposition \ref{proposition : two torus embedding}, there exists a symplectic two torus $T$ in $M$ such that $\langle c_1(TM), [T] \rangle = \lambda \cdot \langle [\omega], [T] \rangle = 0$, which leads to a contradiction. Consequently, we have $c_1(M,\omega) = 0$, i.e., $\lambda = 0$.  
By the same reason as before, if $\lambda$ is negative, then there is no circle subgroup of $G$. In other words, $G$ must be trivial by the effectiveness of the $G$-action.
\end{proof}

We make the following two remarks. First, the reason why we only prove Proposition \ref{proposition : two torus embedding} for the rational case is that our proof relies on the existence of a \textit{generalized moment map} \cite{McD} which can be defined only when $[\omega]$ is rational. If $[\omega]$ is not rational, since the non-degeneracy of $\omega$ is an open condition, we can always perturb a given $\omega$ slightly to another symplectic form $\omega'$ such that $\omega'$ is $G$-invariant and $[\omega']$ is rational. Hence if we apply Proposition \ref{proposition : two torus embedding} to $(M,\omega')$, then there exists a symplectic embedding $i : T^2 \hookrightarrow M$ with respect to the new symplectic structure $\omega'$. But there is no guarantee that the $T^2$-embedding $i$ is symplectic with respect to $\omega$. The authors could not fill the gap of the proof of Proposition \ref{proposition : two torus embedding} in the case when $\omega$ is not rational.

Secondly, unlike Proposition \ref{proposition : two torus embedding} which gives a necessary and sufficient condition for the existence of non-Hamiltonian $G$-action, Proposition \ref{proposition : embedding sphere Hamiltonian case} gives only the sufficient condition for the action to be Hamiltonian. The authors do not know whether the converse of Proposition \ref{proposition : embedding sphere Hamiltonian case} holds or not.

The organization of this paper is as follows : we give an introduction to the theory of Lie group actions on symplectic manifolds in Section 2, and we give the complete proof of Proposition \ref{proposition : two torus embedding} and \ref{proposition : embedding sphere Hamiltonian case} in Section 3.

\section{Symplectic circle Actions}

In this section, we give a brief introduction to symplectic
circle actions. Most of this section is contained in
\cite{Au} or \cite{McS}, but we give a complete proof for
readers who are not familiar with symplectic geometry. Let
$M$ be a $2n$-dimensional smooth closed manifold. A
differential 2-form $\omega$ is called a
\textit{symplectic form} if $\omega$ is closed and
non-degenerate, i.e.,
\begin{itemize}
    \item $d\omega = 0$, and
    \item $\omega^n$ is nowhere vanishing.
\end{itemize}
Let us assume that $G$ is a compact connected Lie group acting effectively on $M$.
A $G$-action on $(M,\omega)$ is called \textit{symplectic} if $\mathcal{L}_{\underline{X}}\omega = 0$ for every $X \in T_eG$ where $\underline{X}$ is the fundamental vector field generated by $X$, i.e., $G$ preserves a symplectic form $\omega$. Equivalently, $G$-action is symplectic if and only if $i_{\underline{X}}\omega$ is closed. In particular, a $G$-action is called \textit{Hamiltonian} if $i_{\underline{X}}\omega$ is exact for every $X \in T_eG$.

Now suppose that the unit circle group $S^1$ acts on
$(M,\omega)$ symplectically, and let $X$ be a fixed generator of
$T_e S^1 \cong \R$. If the action is
Hamiltonian, then there exists a smooth function $H : M
\rightarrow \R$ such that $$ i_{\underline{X}}\omega = dH,
$$ and we call $H$ a \textit{moment map} for the $S^1$-action.
Note that $\mathcal{L}_X H =
\omega(\underline{X},\underline{X}) = 0$, i.e., a moment map
$H$ is $S^1$-invariant.

If the $S^1$-action is symplectic but non-Hamiltonian, then
a moment map does not exist. Nevertheless, there exists an
$\R/\Z$-valued function $\mu : M \rightarrow \R/\Z$ which
locally looks like a moment map when $[\omega] \in
H^2(M;\R)$ is an integral class. We use the notation
$\R/\Z$ instead of $S^1$ to avoid confusion with the acting
group $S^1.$

\begin{definition}\label{definition : generalized moment map}\cite{McD}
Let $(M,\omega)$ be a smooth closed symplectic manifold
such that $\omega$ represents an integral cohomology class
in $H^2(M;\Z)$. Suppose that there is a symplectic
non-Hamiltonian $S^1$-action on $(M,\omega)$. Fix a point
$x_0 \in M$ and define an $\R/\Z$-valued map $\mu : M
\rightarrow \R/\Z$ such that $$\mu(x) := \int_{\gamma_x}
i_{\underline{X}}\omega ~\mod \Z$$
    where $\gamma_x$ is any path $\gamma_x : [0,1] \rightarrow M$ such that $\gamma_x(0) = x_0$ and $\gamma_x(1) = x$.
    We call $\mu$ an \textit{$\R/\Z$-valued moment map} (or a \textit{generalized moment map}) for the action.
\end{definition}

By a direct computation, we can easily check that $[i_{\underline{X}}\omega]$ is an integral class in $H^1(M;\Z)$ so that a generalized moment map given in Definition \ref{definition : generalized moment map} is well-defined. Note that $\mu$ depends on the choice of a base point $x_0$ as in Definition \ref{definition : generalized moment
map}. Consider two distinct points $p$ and $q$ on $M$, and
let $\mu_p$, respectively $\mu_q$, be the $\R/\Z$-valued moment map
with base point $p$, respectively $q$. For any point $x
\in M$, let $\gamma_q^p$ be a path from $q$ to $p$ and
$\gamma_p^x$ be a path from $p$ to $x$ respectively. Then
$$\mu_p(x) - \mu_q(x) = \int_{\gamma_p^x}
i_{\underline{X}}\omega - \int_{\gamma_p^x \circ
\gamma_q^p} i_{\underline{X}}\omega = -\int_{\gamma_q^p}
i_{\underline{X}}\omega = -\mu_q(p) ~\mod \Z.$$ In other
words, $\mu$ is unique up to a constant in $\R / \Z \cong
S^1$. In particular, $d\mu$ is independent of the choice of
a base point. Since $d\mu : TM \rightarrow TS^1
\cong S^1 \times \R$, we may regard $d\mu$ as a
differential 1-form on $M$.

\begin{proposition}\cite{McD}\label{proposition : dmu = ixomega}
Let $\mu : M \rightarrow \R/\Z$ be an $\R/\Z$-valued moment
map for a symplectic non-Hamiltonian circle action on $(M,\omega)$.
Then $\mu$ satisfies $$d\mu = i_{\underline{X}}\omega.$$
\end{proposition}
\begin{proof}
    Let $x \in M$ be any point and let $\mathcal{U}$ be a contractible open neighborhood of $x$. Since $i_{\underline{X}}\omega$ is closed, it is locally exact by Poincar\'{e} lemma so that there exists a smooth function $f : \mathcal{U} \rightarrow \R$ such that $i_{\underline{X}}\omega = df$ on $\mathcal{U}$. Let $\mu$ be an $\R/\Z$-valued moment map with a base point $x_0 \in \mathcal{U}$. Let $\gamma_x$ be a path from $x_0$ to $x$ lying on $\mathcal{U}$. Then $$\mu(x) = \int_{\gamma_x} i_{\underline{X}}\omega = \int_{\gamma_x} df = f(x) - f(x_0) \mod\Z$$ so that $d\mu(x) = df(x) = i_{\underline{X}(x)}\omega_x$ for all $x \in \mathcal{U}$. Since $x$ is chosen arbitrarily, we can conclude that $d\mu = i_{\underline{X}}\omega$ on $M$.
\end{proof}

It is an immediate consequence of Proposition \ref{proposition : dmu = ixomega} that $\mu$ is $S^1$-invariant, since $\mathcal{L}_{\underline{X}} \mu = i_{\underline{X}}d\mu = \omega(\underline{X}, \underline{X}) = 0$.

Now, let us consider a critical point of a moment map $H$,
i.e., $dH(x) = i_{\underline{X}(x)}\omega_x = 0$. Since
$\omega$ is non-degenerate on $M$, $x$ is a critical point
of $H$ if and only if $\underline{X}(x) = 0$, i.e., $x$ is a
fixed point of given $S^1$-action. It is also true for a
non-Hamiltonian case, i.e., $x$ is a critical point of an
$\R/\Z$-valued moment map $\mu$ if and only if $x$ is a
fixed point of the action by Proposition \ref{proposition :
dmu = ixomega}. Hence we have the following.

\begin{proposition}\label{proposition : crit = fixed points}\cite{Au}
    Let $(M,\omega)$ be a smooth closed symplectic manifold equipped with a symplectic, respectively Hamiltonian, circle action. If $\mu$, respectively $H$, is an $\R/\Z$-valued moment map, respectively moment map, of given action, then $x \in M$ is a critical point of $\mu$, respectively $H$, if and only if $x$ is a fixed point of the action.
\end{proposition}

One of the most important property of symplectic geometry is that, for every point $p \in M$, there exists a neighborhood $\mathcal{U}_p$ of $p$ with a local coordinate system $(\R^{2n}, x_1, y_1, \cdots, x_n, y_n)$ such that $\omega|_{\mathcal{U}_p} = \sum dx_i \wedge dy_i$, i.e., a local symplectic structure for each point is isomorphic to the standard symplectic structure of $\R^{2n}$ so that symplectic geometry is locally the same as a linear symplectic geometry on $\R^{2n}$ with the standard symplectic structure $\sum dx_i \wedge dy_i$. This is known as the \textit{Darboux theorem}. Similarly, there is an equivariant version of the Darboux theorem as follows.

\begin{theorem}[Equivariant Darboux theorem]\label{theorem : equivariant Darboux}
    Let $(M,\omega)$ be a symplectic manifold and let $G$ be a compact Lie group.
    Suppose that there is a symplectic $G$-action on $(M,\omega)$. For each fixed point $p$, there exists a neighborhood $\mathcal{U}_p$ together with a local coordinate system $(x_1,y_1, \cdots, x_n, y_n)$ such that
    \begin{itemize}
        \item $\omega|_{\mathcal{U}_p} = \frac{1}{2i} \sum dz_j \wedge d\bar{z_j}$ with $z_j = x_j + iy_j$, and
        \item $G$-action is linear with respect to $(z_1, \cdots, z_n)$. In particular if $G = S^1$, then there is a sequence of integers $\lambda_1, \cdots, \lambda_n$ such that the action is expressed as
        $$t \cdot (z_1, \cdots, z_n) = (t^{\lambda_1}z_1, \cdots, t^{\lambda_n}z_n)$$ for every $t \in S^1$.
    \end{itemize}
\end{theorem}

Now, let $p$ be a fixed point of symplectic $S^1$-action on
$2n$-dimensional symplectic manifold $(M,\omega)$. By the
equivariant Darboux theorem, there exists a local
coordinate system $(\mathcal{U}_p, z_1,\cdots,z_n)$
centered at $p$ and a sequence of integers $\lambda_1,
\cdots, \lambda_n$ such that the $S^1$-action is expressed
as
\[
t \cdot (z_1, \cdots, z_n) = (t^{\lambda_1}z_1, \cdots,
t^{\lambda_n}z_n)
\]
for every $t \in S^1$. By solving $i_{\underline{X}}\omega
= dH$ on $\mathcal{U}_p$ with $\omega = \frac{1}{2i} \sum
dz_j \wedge d\bar{z_j}$, we get

\begin{equation} \label{equation: local form of H}
H(z_1,\cdots, z_n) = \text{constant}+ \frac{1}{2} \sum
\lambda_j |z_j|^2.
\end{equation}
Therefore, we have the following corollary.

\begin{corollary}\cite{Au}\label{corollary : collary of equiv Darboux}
    Let $H : M \rightarrow \R$ be a moment map on $(M,\omega)$. Then $H$ is a Morse-Bott function.
    Similarly, if $\mu : M \rightarrow \R/\Z$ is an $\R/\Z$-valued moment map, then $\mu$ is an $\R/\Z$-valued Morse-Bott function. In either case, a Morse index of any critical submanifold of $H$ or $\mu$ is even.
\end{corollary}

\begin{proof}
We need to show two things : (1) the critical point set is
an embedded submanifold of $(M,\omega)$, and (2) the
Hessian of $H$ at $p$ is non-degenerate along the normal
direction of a critical submanifold containing $p$. The
first claim is obvious since a sub-coordinate system
$$\{(z_1, \cdots, z_n) \in \mathcal{U}_p ~|~ z_j = 0
~\mathrm{if} ~\lambda_j \ne 0\}$$ gives a coordinate system
of a critical submanifold near $p$. The Hessian of $H$ at
$p$ is a diagonal matrix is given by
    \begin{displaymath}
        \left(\begin{array}{cccccc}
        \lambda_1 & 0 & 0 & \cdots & 0 & 0\\

        0 & \lambda_1 & 0 & \cdots & 0 & 0\\

        \vdots & \vdots & \vdots & \vdots & \vdots & \vdots\\

        \vdots & \vdots & \vdots & \vdots & \vdots & \vdots\\

        0 & 0 & \cdots & 0 & \lambda_n & 0 \\

        0 & 0 & \cdots & 0 & 0 & \lambda_n \\

    \end{array}\right)
    \end{displaymath}
    so that it finishes the proof.
\end{proof}

Now, recall some basic Morse-Bott theory as follows. Let $f
: M \rightarrow \R$ be a Morse-Bott function on a closed
manifold $M$, and let $M_t = \{ p \in M |~ f(p) \leq t \}$ for every $t \in \R$. Suppose that $a$ and $b$ are regular values
of $f$ such that there exists a unique critical value $c$
with $a < c < b$. Let $C_1, \cdots , C_r$ be connected
components of the critical submanifold lying on
$H^{-1}(c)$. According to classical Morse-Bott theory,
$M_b$ is homotopy equivalent to
$$M_a \cup_{\phi_1} D(\nu^-(C_1)) \cup_{\phi_2}
D(\nu^-(C_2)) \cdots \cup_{\phi_r} D(\nu^-(C_r))$$ where
$\nu^-(C_j)$ is a negative normal bundle over $C_j$,
$D(\nu^-(C_j))$ is a disk bundle of $\nu^-(C_j)$, and
$\phi_j$ is an attaching map from a sphere bundle
$S(\nu^-(C_j)) = \partial D(\nu^-(C_j))$ to $H^{-1}(a)$.
Note that each $S(\nu^-(C_j))$ is an $S^{k_j-1}$-bundle
over $C_j$ where $k_j = \mathrm{ind}(C_j)$ is a Morse index
of $C_j$. In particular, $S(\nu^-(C_j))$ is connected if
and only if $k_j \neq 1$.

\begin{proposition}\cite{Au}\label{proposition : connectivity of level set : Hamiltonian case}
Let $H$ be a moment map on a (possibly non-compact)
connected symplectic manifold $(M,\omega)$. Then every
level set of $H$ is empty or connected.
\end{proposition}

\begin{proof}
For the sake of simplicity, let $M(a,b) := H^{-1}((a,b))$ for $a,b \in \R$. 
Let us choose any $S^1$-invariant metric $\langle \cdot, \cdot \rangle$ on $M$ so that the gradient flow $\bigtriangledown H $ of $H$ is defined as $$ dH(X) = \langle \bigtriangledown H, X \rangle $$ for every smooth vector field $X$ on $M$. Suppose that there exists a regular value $r \in \R$ such
that $H^{-1}(r)$ is non-empty and disconnected. By the
connectivity of $M$, there exists a smallest $s \in \R^+$
such that $M(r-s,r+s)$ is connected. Note that $r+s$ or
$r-s$ is a critical value of $H$, otherwise $M(r-s,r+s)$ is
diffeomorphic to $M(r-(s-\epsilon), r+(s-\epsilon))$ along
the gradient flow of $H$ for a sufficiently small $\epsilon
>0$ so that it contradicts our assumption ``smallest $s$''.

    Without loss of generality, we may assume that $c = r+s$ is a critical value of $H$, Let $C_1, \cdots , C_r$ be connected components of the critical submanifold lying on $H^{-1}(c)$.
    Then there exists some $C_j$ such that $D(\nu^-(C_j))$ connects two disconnected components of $H^{-1}(c-\epsilon)$ via the attaching map $\phi_j : S(\nu^-(C_j)) \rightarrow H^{-1}(c-\epsilon)$, i.e., the index of $C_j$ should equal to one. Since every critical submanifold of $H$ has even index by Corollary \ref{corollary : collary of equiv Darboux}, such $C_j$ does not exist.

    Similarly, if $c = r-s$ is a critical value of $H$, then there exists some $C_j$ of co-index one, but there is no such $C_j$ by Corollary \ref{corollary : collary of equiv Darboux}. Hence it completes the proof.

\end{proof}

\begin{proposition}\cite{Au}\label{proposition : same number of connected components}
Let $(M,\omega)$ be a $2n$-dimensional closed connected
symplectic manifold equipped with a symplectic
non-Hamiltonian circle action. Suppose $[\omega]$ is
an integral class in $H^2(M;\Z)$, and let $\mu : M
\rightarrow \R/\Z$ be an $\R/\Z$-valued moment map defined
in Definition \ref{definition : generalized moment map}.
Then there is no critical submanifold of index zero nor co-index zero.
In particular, every level set is non-empty and the number of
connected components of $\mu^{-1}(t)$ is constant for all $t
\in S^1$.
\end{proposition}

\begin{proof}
Let $r \in \R / \Z$ be a regular value of $\mu : M \rightarrow
\R/\Z$ and let $N = \mu^{-1}(\R / \Z / \{r\})$ which is an open
subset in $M$. With the induced $S^1$-action on $N$, we may
regard $N$ as a non-compact Hamiltonian $S^1$-manifold with
a moment map $H = \mu|_N : N \rightarrow \R/\Z - \{r\}
\cong (0,1)$. Let $N_1, N_2, \cdots, N_k$ be connected
components of $N$ and we denote by $H_j : N_j \rightarrow
(0,1)$ the restriction of $H$ onto $N_j$ so that $H_j$ is a
moment map on $(N_j, \omega|_{N_j})$. By Proposition
\ref{proposition : connectivity of level set : Hamiltonian
case}, every level set of $H_j$ is empty or connected for
every $j=1,2,\cdots, k$.

    Firstly, we claim that each $H_j$ is surjective. If not, $H(N_j)$ is either a half-closed interval of the form $[s,1)$ or $(0,s]$ for some $s \in (0,1)$ or a closed interval $[a,b] \subset (0,1)$. If $H(N_j) = [a,b]$, then $N_j$ is also a connected component of $H^{-1}([a,b])$ so that $N_j$ is both open and closed itself in $M$ so that $N_j = M$ by the connectivity of $M$, which contradicts to the assumption that the action is non-Hamiltonian.
    If $H(N_{j_1}) = [s_1,1)$ for some $j_1 \in \{1,2,\cdots, k \}$, then let us consider the closure $\overline{N}_{j_1}$ whose boundary $\partial \overline{N}_{j_1}$ is some connected component, namely $B_1$, of $\mu^{-1}(r)$. Since $r$ is regular, there is a connected component $N_{j_2}$ of $N$ such that $H_{j_2}^{-1}(t)$ attains $B_1$ as $t \rightarrow 0$. If the image of $N_{j_2}$ for $H_{j_2}$ is a half-closed interval of the form $(0,s_2]$ for some $s_2 \in (0,1)$, then $\overline{N}_{j_1} \cup \overline{N}_{j_2}$ is connected and both open and closed in $M$ so that we get $\overline{N}_{j_1} \cup \overline{N}_{j_2} = M$. Then the moment map $\mu$ factors such that
    $$ \mu : \overline{N}_{j_1} \cup \overline{N}_{j_2} = M \stackrel{H_{j_1,j_2}} \longrightarrow [s_1, 1+s_2] \stackrel{ / \Z} \longrightarrow S^1 $$
    where $H_{j_1,j_2}$ maps $x \in N_{j_1}$ to $H_{j_1}(x)$, $y \in \mu^{-1}(r)$ to $1$, and $z \in N_{j_2}$ to $1 + H_{j_2}(z)$. Then $i_{\underline{X}}\omega = d(/ \Z) \circ dH_{j_1,j_2}$, but $d(/ \Z)$ is the identity map so that $dH_{j_1,j_2} = i_{\underline{X}}\omega$, i.e., the action is Hamiltonian which contradicts to our assumption. Hence $H(N_{j_2}) = (0,1)$, i.e., $H_{j_2}$ is surjective so that $H_{j_2}^{-1}(t)$ attains some connected component $B_2 \neq B_1$ of $\mu^{-1}(r)$ as $t \rightarrow 1$. Take $N_{j_3}$ such that $H_{j_3}^{-1}(t)$ attains $B_2$ as $t \rightarrow 0$. Then we can easily show that $H_{j_3}$ is surjective by a similar reason. Hence we get a infinite sequence of pairwise distinct connected components $B_1, B_2, \cdots$, but it contradicts that the number of connected components of $\mu^{-1}(r)$ is finite by the compactness of $M$. Therefore, $H_j$ is surjective for every $j$. In particular, there is no critical submanifold of index 0 nor co-index 0.

    To complete the proof, recall that when a level set of $\mu$ passes through some critical value $c \in S^1$ such that $\mu^{-1}(c)$ does not contain a critical submanifold of index 0, 1,  co-index 0, nor co-index 1, then the number of connected components doesn't change, i.e., $\mu^{-1}(c+\epsilon)$ and $\mu^{-1}(c-\epsilon)$ have the same number of connected components, see the proof of Proposition \ref{proposition : connectivity of level set : Hamiltonian case}.
    Since an index and co-index of any critical component is even and there is no critical component of index 0 nor co-index 0, every level set has the same number of connected components.
\end{proof}

\section{Proof of the Main theorem}
Let $G$ be a compact Lie group, and suppose 
$(M,\omega)$ is a closed symplectic manifold equipped with
an effective symplectic $G$-action. In this section, we
prove Theorem \ref{theorem : main}. To determine whether a
given $G$-action is Hamiltonian or not, it is enough to
check it for every circle subgroup of $G$ since every
element $g \in G$ is contained in some maximal torus of
$G$. The following proposition characterizes a
non-Hamiltonian circle action in terms of equivariant
symplectic embedding of two-torus.

\begin{proposition}[Proposition \ref{proposition : two torus embedding}]
    Let $(M,\omega)$ be a smooth closed symplectic manifold equipped with a smooth $S^1$-action preserving $\omega$.
    Suppose that $[\omega]$ is a rational class in $H^2(M;\Q)$. Then the action is non-Hamiltonian if and only if there exists an $S^1$-equivariant symplectic
    embedding of 2-torus $i \colon T^2 \hookrightarrow M$, where the circle acts freely on the left factor of $T^2 \cong S^1 \times S^1$.
    In particular, the normal bundle of $i(T^2)$ in $M$ is trivial so that $\langle c_1(M,\omega), [i(T^2)] \rangle = 0$.
\end{proposition}

\begin{proof}

Let $T^2 = S^1 \times S^1$ be a two-torus and assume that
an $S^1$-action on $T^2$ is given by
\[
t \cdot (t_1,t_2) = (t \cdot t_1,t_2)
\]
for any $t \in S^1$ and $(t_1,t_2) \in T^2$. Suppose there exists an
$S^1$-equivariant symplectic embedding $i \colon T^2 \hookrightarrow M$. If the given $S^1$-action on $(M,\omega)$ is Hamiltonian with a moment map $H : M \rightarrow \R$, then it is straightforward that the restriction of the $S^1$-action on $(T^2, i^*\omega)$ is automatically Hamiltonian with a moment map $i^*H = H \circ i : T^2 \rightarrow \R$. But it contradicts to the assumption that $S^1$-action on $T^2$ is free, because there must be at least two fixed points, namely the maximum and the minimum of $i^*H$. Hence the $S^1$-action on $(M,\omega)$ cannot be Hamiltonian.

Conversely, suppose that the given symplectic $S^1$-action on $(M,\omega)$ is non-Hamiltonian.
By our assumption, there exists a natural number $N$ big enough so that $N \cdot \omega$ is integral and we still denote by $\omega$ the new symplectic form $N \cdot \omega$.
Obviously, the given $S^1$-action is symplectic with respect to the new symplectic form $\omega$ so that there exists an $\R/\Z$-valued moment map $\mu : M \rightarrow \R/\Z$
satisfying $i_{\underline{X}}\omega = d\mu$ by Proposition \ref{proposition : dmu = ixomega}.
Let $M_{(1)}$ be the set of all points in $M$ with the trivial isotropy subgroup.
Now, suppose that there exists a smoothly embedded loop $\sigma: S^1 = [0,n] /_{\mathord{0 \sim n}} \longrightarrow M_{(1)}$ for some
$n \in \N$ satisfying the following conditions:
\begin{itemize}
\item[(a)] $\mu \big( \sigma(r) \big) = [r] \in \R / \Z$ for each $r \in [0,n]$,
\item[(b)] $\sigma(0) = \sigma(n)$, i.e., the image of $\sigma$ is a loop, and
\item[(c)] for $r \ne r^\prime$ with $(r,r') \neq (0,n)$, $t \cdot \sigma(r) \ne \sigma(r^\prime)$ for any $t \in S^1$.
\end{itemize}
If such $\sigma$ exists, then we can define a smooth embedding of 2-torus as follows :
\[
i : S^1 \times [0,n] /_{\mathord{0 \sim n}} \cong T^2 \rightarrow M, \quad (t_1, t_2) \longmapsto t_1 \cdot \sigma(t_2).
\]
It is straightforward that $i$ is $S^1$-equivariant for the $S^1$-action on $T^2$ by $t \cdot (t_1, t_2) = (t
\cdot t_1, t_2).$ To show that $i(T^2)$ is a symplectic submanifold with respect to the induced symplectic structure, let us define a smooth vector field $Y$ on $i(T^2)$ as
follows:
\begin{equation} \label{equation: vector field}
Y \Big( i(t_1,t_2) \Big) ~ := ~ \frac{ d(t_1 \cdot \sigma(t))}{dt} \Big|_{t = t_2}.
\end{equation}
Since $\sigma$ has no critical point, it is straightforward that $Y$ has no zero and $\underline{X}(p)$  and
$Y(p)$ span the tangent space $T_p ~ i(T^2)$ for every $p
=i(t_1, t_2) \in i(T^2)$. Also,
\begin{alignat*}{2}
\omega \big( \underline{X}(p), Y(p) \big) &~=~ d\mu \big(
Y(p) \big) & & \qquad \text{ by definition of } \mu \\
&~=~ d\mu \Big( \frac{d(t_1 \cdot \sigma(t))}{dt}
\Big|_{t = t_2} \Big) & & \qquad \text{ by definition of } Y  \\
&~=~ d\mu \Big( \, \frac{d \sigma}{dt} \Big|_{t = t_2}
\Big) & & \qquad \text{ by $S^1$-invariance of } \mu
\\
&~=~ 1 & & \qquad \text{ by } (a).
\end{alignat*}
So, $\omega$ is non-degenerate on $i(T^2)$, i.e., $i$ is a
symplectic embedding. 

Now, we need to prove that such smoothly embedded loop $\sigma$ in $M$ actually exists.

\begin{lemma} \label{lemma : principal orbit}
$M_{(1)}$ is path-connected and open dense in $M.$
\end{lemma}
\begin{proof}
Let $\Z_n$ be the cyclic subgroup of $S^1$ of order $n$, and we denote by $M^{\Z_n}$ the set of all points in $M$ with the isotropy subgroups equal to $\Z_n$.
Since $M$ is compact, it is well known there are at most finitely many $n$'s, say $n_1, n_2, \cdots, n_k$ such that $M^{\Z_n} \neq \emptyset$, see \cite[Proposition IV.1.2]{Br}. Let  $M^{S^1}$ be the set of all points in $M$ fixed by $S^1$.
Since $M^{S^1} \subset M^{\Z_n}$ for every integer $n >1$, we have 
\begin{equation} \label{equation: trivial stabilizer}
M_{(1)} ~=~ M - \bigcup_{n>1} M^{\Z_n} - M^{S^1} ~=~ M -
\bigcup_{n>1} M^{\Z_n}.
\end{equation}
Furthermore, 
$M^{\Z_n}$ and $M^{S^1}$ are closed symplectic submanifolds of $M$ with the
induced symplectic form by the equivariant Darboux theorem
\ref{theorem : equivariant Darboux}. Thus the set
$M^{S^1} \cup (\bigcup_{n>1} M^{\Z_n}) $ is the union of
closed submanifolds with codimensions at least two, in particular
$M_{(1)}$ is path-connected and open dense in $M$.
\end{proof}

\begin{corollary} \label{corollary : level set open dense}
For any regular value $t_0 \in \R/\Z$ of $\mu,$ the subset
$M_{(1)} \cap \mu^{-1}(t_0)$ is open dense in $\mu^{-1}(t_0)$.
\end{corollary}

\begin{proof}
The openness is obvious since $M_{(1)}$ is open.
Also, it is not hard to show that $M^{\Z_n} \cap \mu^{-1}(t_0)$ is of codimension at least two in $\mu^{-1}(t_0)$ for each $n > 1$.
Since $t_0$ is chosen to be regular, we have $\mu^{-1}(t_0) \cap M^{S^1}$ is empty so that $M_{(1)} \cap \mu^{-1}(t_0) = \mu^{-1}(t_0) - \bigcup_{n>1} M^{\Z_n}$ is a dense subset of $\mu^{-1}(t_0)$.

\end{proof}

\begin{remark}
The set $M_{(1)} \cap \mu^{-1}(t_0)$ is not necessarily
path-connected in Corollary \ref{corollary : level set open
dense}.
\end{remark}

Without loss of generality, we assume that $0 \in \R/\Z$ is
a regular value of $\mu.$ For our convenience, we use the following terminology : for a fixed $S^1$-invariant metric $h$, we say that a smooth path $\gamma : [a,b] \rightarrow M_{(1)}$ \textit{winds} along $\R / \Z$ if $h(\bigtriangledown \mu(t), \big(\mu \circ \gamma \big)^\prime (t)) >0$ for every $t \in [a,b]$, which means that the vector field generated by $\gamma$ is a \textit{gradient-like vector field} of $\mu$ with respect to $h$. We call such a path $\gamma$ a \textit{winding path}.
Also, we say that two points $x, y \in M_{(1)}$ are \textit{winding path-connected} if there exists a path
$\gamma: [a,b] \rightarrow M^{(1)}$ from $\gamma(a)=x$ to $\gamma(b)=y$ which winds along $\R / \Z$.

\begin{lemma} \label{lemma : winding-path connected}
For $x \in M_{(1)} \cap \mu^{-1}(0)$, there exists some $y \in M_{(1)} \cap \mu^{-1}(0)$ such that $x$ and $y$ are winding path-connected.
\end{lemma}

\begin{proof}
 Let $J$ be an $S^1$-invariant almost complex structure compatible\footnote{We say that an almost complex structure $J$ on $(M, \omega)$ is \textit{compatible with $\omega$} if  (1) $\omega(\cdot, \cdot) = \omega(J\cdot, J\cdot)$  and (2) $\omega(\cdot, J\cdot)$ is a Riemannian metric. Such $J$ always exists, see \cite{McS} for the details.} with $\omega$.
We denote by $\langle \cdot, \cdot \rangle = \omega(\cdot, J\cdot)$ the induced $S^1$-invariant metric on $M$ so that the gradient vector field $\bigtriangledown \mu$ on $M$ with respect to $\langle \cdot, \cdot \rangle$ is defined as $d\mu = \langle \bigtriangledown \mu, \cdot \rangle$.
Since $d\mu = i_{\underline{X}}\omega = \langle \bigtriangledown \mu, \cdot \rangle$, we can easily see that $\bigtriangledown \mu = J\underline{X}$.
Note that $\bigtriangledown \mu$ commutes with the $S^1$-action since $J$, $\langle \cdot, \cdot \rangle$, and $\mu$ are chosen to be $S^1$-invariant. Hence the one-parameter group action generated by $\bigtriangledown \mu$ preserves their isotropy subgroups, which means that if $x \in M$ is fixed by some subgroup $H \subset S^1$, then any point $y$ in the orbit $\{ (\exp{t {\bigtriangledown \mu}}) \cdot x \}_{t \in \R}$ is fixed by $H$. In particular, the one-parameter group action generated by $\bigtriangledown \mu$ acts on $M_{(1)}$.

Now, fix $x \in M_{(1)}$ and let us consider the integral curve along $\bigtriangledown \mu$ passing through $x$
\begin{displaymath}
    \begin{array}{ccccc}
        \gamma_x & : & \R & \longrightarrow & M \\
                   &   &  t & \mapsto        & (\exp{t {\bigtriangledown \mu}}) \cdot x. \\
    \end{array}
\end{displaymath}
If $\gamma_x$ is a winding path from $x$ to some point $y \in M_{(1)} \cap \mu^{-1}(0)$, then there is nothing to prove.
If not, then $\lim_{t \rightarrow \infty} \mu(\exp{t {\bigtriangledown \mu}} \cdot x) = t_0 $ for some $t_0 \in S^1$, which is equivalent to saying that $$\lim_{t \rightarrow \infty} \gamma_x(t) = \lim_{t \rightarrow \infty} (\exp{t {\bigtriangledown \mu}}) \cdot x = p $$ for some fixed point $p \in M^{S^1}$.
By the equivariant Darboux theorem \ref{theorem : equivariant Darboux}, there exists a local coordinate system $(\mathcal{U}_p, z_1,\cdots,z_n)$ centered at $p$ such that

 $$\omega|_{\mathcal{U}_p} = \frac{1}{2i} \sum dz_j \wedge d\bar{z_j},$$ and
 $$t \cdot (z_1, \cdots, z_n) = (t^{\lambda_1}z_1, \cdots, t^{\lambda_n}z_n)$$ for every $t \in S^1$ where $\lambda_1, \cdots, \lambda_n$ are weights of the tangential $S^1$-representation  on $T_pM$.
Let $C_p$ be the fixed connected component containing $p$ and let $\nu_+$ and $\nu_-$ be subsets of $\mathcal{U}_p$ given by
\begin{itemize}
    \item $\nu_+ = \{(z_1, \cdots, z_n) \in \mathcal{U}_p ~|~ z_j = 0 \text{ if } \lambda_j \le 0\}$, and
    \item $\nu_- = \{(z_1, \cdots, z_n) \in \mathcal{U}_p ~|~ z_j = 0 \text{ if } \lambda_j \ge 0\}. $
\end{itemize}
Note that $\{(z_1, \cdots, z_n) \in \mathcal{U}_p ~|~ z_j \neq 0 \text{ if } \lambda_j \neq 0\} \subset M_{(1)}$ since the $S^1$-action is effective.
Hence we may perturb $\gamma_x$ smoothly in a sufficiently small neighborhood $\mathcal{U}$ of $C_p$ to get a new gradient-like flow $\widetilde{\gamma}_x$ on $M_{(1)}$ connecting the gradient flow $\gamma_x$ to $\gamma_q$ for some $q \in M_{(1)}$ as in Figure \ref{figure: gradient flow} below :
\begin{figure}[ht]
\begin{center}
\begin{pspicture}(-3,-2)(3,2) \footnotesize

\psline(-2,0)(2,0) \psline(0,-2)(0,2)

\psline[linewidth=1.2pt,arrowsize=8pt]{->}(2,0)(1.4,0)
\psline[linewidth=1.2pt](2,0)(1.2,0)
\psline[linewidth=1.2pt](0.2,1)(0.1,2)
\psline[linewidth=1.2pt,arrowsize=8pt]{->}(0.2,1)(0.15,1.5)
\psline[linewidth=1pt,arrowsize=4pt]{->}(0.2,0)(0.15,0)

\psdot(0,0)
\psarc[linewidth=1.2pt,showpoints=false](1.2,1){1}{180}{270}
\psdot(0.15,1.5)
\pscircle[linestyle=dotted](0,0){1.51}

\uput[r](2,0){$\nu_-$} \uput[u](0,2){$\nu_+$}
\uput[d](1.8,0.15){$\downarrow$}
\uput[d](1.8,-0.2){$\gamma_x$} \uput[dl](0,0){$C_p$}
\uput[r](0.2,0.6){$\widetilde{\gamma}_x \subset M_{(1)}$}
\uput[r](0.15,1.5){$q$}
\uput[r](0.03,1.8){$\rightsquigarrow \gamma_q = \exp(t {\bigtriangledown \mu}) \cdot q$}
\uput[r](1.1,-1.1){$\mathcal{U}$}
\end{pspicture}
\end{center}
\caption{ \label{figure: gradient flow} Perturbing $\gamma_x$ to get $\widetilde{\gamma}_x$}
\end{figure}
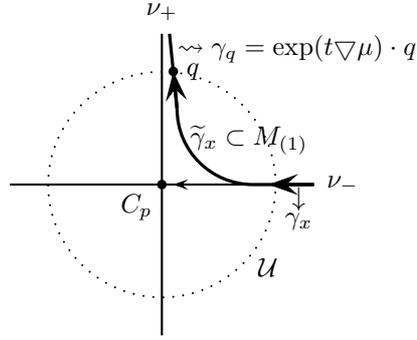
If $\widetilde{\gamma}_x$ is a winding path from $x$ to some point $y \in M_{(1)} \cap \mu^{-1}(0)$, then it is done. If not, we apply this procedure finitely many times so that
we can get a winding path connecting $x$ to $M_{(1)} \cap \mu^{-1}(0)$ as is desired.
\end{proof}

\begin{lemma} \label{lemma : whole component}
If a point $x$ in $M_{(1)} \cap \mu^{-1}(0)$ is winding path-connected to some point $y$ in a connected component $C$ of $M_{(1)} \cap \mu^{-1}(0)$,
then $x$ is winding path-connected to every point in $C$.
\end{lemma}

\begin{proof}
For $x \in M_{(1)} \cap \mu^{-1}(0)$, let $\gamma_x$ be a winding path connecting $x$ to $y \in M_{(1)} \cap \mu^{-1}(0)$.
Let $C_y$ be the connected component of $M_{(1)} \cap \mu^{-1}(0)$ containing $y$.
Fix $z \in C_y$ and let $\sigma : [0,1] \rightarrow M_{(1)} \cap \mu^{-1}(0)$ be a smooth path with $\sigma(0) = y$ and $\sigma(1) = z$.
Since $0$ is chosen to be a regular value of $\mu$, the gradient vector field $\bigtriangledown \mu$ is non-zero on $\mu^{-1}(0)$. Also, we can find a sufficiently small $\epsilon$ such that any $r \in [-\epsilon,0] \subset S^1$ is a regular value of $\mu$.

\begin{figure}[ht]
\begin{center}
\begin{pspicture}(-5,-2.0)(5,2.2) \footnotesize

\psellipse(0,1)(2,0.5)
\psellipse[linestyle=dotted](0,-1)(2,0.5)
\pscurve[linewidth=0.5pt, linecolor=blue] (-1.5,1)(-1.0,1.2)(-0.5,1)(0.5,0.8)(1.5,1)
\pscurve[linewidth=0.5pt, linecolor=red] (-1.5,-1)(-1.0,-0.8)(-0.5,-1)(0.5,-1.2)(1.5,-1)

\pscurve[linewidth=0.3pt,arrowsize=2pt]{->}(0.5,-1.2)(0.6,-1.6)(0.7,-1.9) \uput[d](0.7,-1.9){$\widetilde{\sigma}$}
\psdots(-1.5,1)(-1.5,-1)(1.5,1)(1.5,-1)
\pscurve[linewidth=0.7pt, arrowsize=4pt]{->}(1.7,-1.5)(1.5,-1)(1.5,1)(1.7,1.5)
\pscurve[linewidth=0.7pt, arrowsize=4pt]{->}(-1.7,-1.5)(-1.5,-1)(-1.6,-0.5)(-1.4,0)(-1.5,1)(-1.7,1.5)
\pscurve[linewidth=0.1pt, linestyle=dashed,arrowsize=4pt]{->}(1.5,-1.5)(1.3,-1)(1.3,1)(1.5,1.5)
\pscurve[linewidth=0.1pt, linestyle=dashed,arrowsize=4pt]{->}(1.1,-1.5)(0.9,-1)(0.9,1)(1.1,1.5)
\pscurve[linewidth=0.1pt, linestyle=dashed,arrowsize=4pt]{->}(0.7,-1.5)(0.6,-1)(0.6,1)(0.7,1.5)
\pscurve[linewidth=0.1pt, linestyle=dashed,arrowsize=4pt]{->}(0.2,-1.5)(0.2,-1)(0.2,1)(0.2,1.5)
\pscurve[linewidth=0.1pt, linestyle=dashed,arrowsize=4pt]{->}(-1.5,-1.5)(-1.3,-1)(-1.4,-0.5)(-1.2,0)(-1.3,1)(-1.5,1.5)
\pscurve[linewidth=0.1pt, linestyle=dashed,arrowsize=4pt]{->}(-1.1,-1.5)(-0.95,-1)(-1.0,-0.5)(-0.85,0)(-0.95,1)(-1.1,1.5)
\pscurve[linewidth=0.1pt, linestyle=dashed,arrowsize=4pt]{->}(-0.8,-1.5)(-0.68,-1)(-0.7,-0.5)(-0.6,0)(-0.65,1)(-0.8,1.5)
\pscurve[linewidth=0.1pt, linestyle=dashed,arrowsize=4pt]{->}(-0.3,-1.5)(-0.19,-1)(-0.2,-0.5)(-0.15,0)(-0.15,1)(-0.3,1.5)
\psline[linewidth=0.3pt,arrowsize=2pt]{->}(1.4,0)(1.8,0) \uput[r](1.7,0){$\gamma_x$}
\psline[linewidth=0.3pt,arrowsize=2pt]{->}(0.3,0.8)(0.5,1.45)(0.7,2)(0.8,2.1) \uput[r](0.7,2.1){$\sigma$}

\pscurve[linewidth=0.5pt, linecolor=gray](-1.5,1)(-1.4,0.7)(1.4,-0.7)(1.5,-1)
\psline[linewidth=0.3pt,arrowsize=2pt]{->}(1.4,-0.7)(1.8,-0.6) \uput[r](1.7,-0.6){$\widetilde{\gamma}_x$}
\uput[r](1.5,1){$y$}
\uput[r](1.5,-1){$y'$}
\uput[l](-1.45,1){$z$}
\uput[l](-1.4,-1){$z'$}
\uput[l](-2,1.2){$M_{(1)} \cap \mu^{-1}(0)$}
\uput[l](-2,-0.8){$M_{(1)} \cap \mu^{-1}(-\epsilon)$}

\end{pspicture}
\end{center}
\caption{ \label{figure: every point winding path-connected}}
\end{figure}
Let $\widetilde{\sigma}$ be the intersection of $M_{(1)} \cap \mu^{-1}(-\epsilon)$ and the trajectory of $\sigma$ under the infinitesimal action generated by $\exp (\bigtriangledown \mu)$.
Let $y'$ be the preimage of $y$ for the infinitesimal action in $M_{(1)} \cap \mu^{-1}(-\epsilon)$. Then we can easily see that there exists a homotopy in $M_{(1)} \cap \mu^{-1}([-\epsilon,0])$ from $\gamma_x$ to a winding path $\widetilde{\gamma}_x$ connecting $y'$ and $z$ (see Figure \ref{figure: every point winding path-connected}).
Consequently, $x$ is winding path-connected to $z$.

\end{proof}

To complete the proof of Proposition \ref{proposition : two torus embedding}, pick a point $x_1$ in a connected component $C_1$ of $M_{(1)} \cap \mu^{-1}(0).$
By Lemma \ref{lemma : winding-path connected}, $x_1$ is winding path-connected to some point $x_2$ in some connected component $C_2$ of $M_{(1)} \cap \mu^{-1}(0)$.
If $C_1=C_2$, then $x_1$ is winding path-connected to $x_1$ itself by Lemma \ref{lemma : whole component} and it satisfies the conditions (a),(b), and (c) automatically.
If $C_1 \neq C_2$, then we can find a point $x_3$ in some connected component $C_3$ of $M_{(1)} \cap \mu^{-1}(0)$ which is winding path-connected to $x_2$.
If $C_1 = C_3$, respectively $C_2 = C_3$, then we can take $x_3 = x_1$, respectively $x_3 = x_2$, so that we can obtain a smooth winding loop satisfying (a),(b), and (c).
Applying the above process inductively, we obtain a smooth winding loop satisfying (a),(b), and (c). Then the following lemma finishes the proof.
\end{proof}

\begin{lemma} \label{lemma : bundle over torus is trivial}
Let $E$ be an $S^1$-equivariant complex vector bundle of rank $k$ over $T^2$ such that the induced $S^1$-action on the zero-section is free. Then $E$ is trivial.
\end{lemma}

\begin{proof}

Let $Z \cong T^2$ be the zero section of $E$. Note that if the induced $S^1$-action on the zero section $Z$ is free, then the given $S^1$-action on the total space $E$ is free.
Now, let us consider a following diagram.
$$\xymatrix{E \ar[r]^{q} \ar[d]^{\pi}& E/{S^1} \ar[d]^{\pi'}\\ Z \ar[r]_{q^\prime}& Z/{S^1}}$$
Since the action is free, $\pi^\prime \colon E/{S^1}
\rightarrow Z/{S^1} \cong S^1$ is a complex vector bundle
of rank $k$ over $S^1$ and the quotient map $q$ is a bundle
morphism. Note that any complex vector bundle over $S^1$ is
trivial, since the structure group $U(n)$ is connected.
Therefore $E/{S^1} \cong S^1 \times \C^k$ and hence $E
\cong q^{\prime *}(E/{S^1}) \cong Z \times \C^k$.
\end{proof}

Now, let us consider the Hamiltonian case. The following proposition characterizes a Hamiltonian circle action in terms of equivariant symplectic embedding of two-spheres.

\begin{proposition}[Proposition \ref{proposition : embedding sphere Hamiltonian case}]
    Suppose $(M,\omega)$ is a smooth closed symplectic manifold with a Hamiltonian circle action.
    Then there exist an $S^1$-equivariant symplectic embedding of 2-sphere $S$ satisfying $\langle c_1(TM), [S] \rangle > 0$.
\end{proposition}

\begin{proof}
    Let $H : M \rightarrow \R$ be a moment map of the circle action. Let $Z_{\min}$, respectively $Z_{\max}$, be a critical submanifold which attains the minimum, respectively maximum, of $H$. Let $g$ be an $S^1$-invariant Riemannian metric defined by $g(X,Y) = \omega(X, JY)$ where $J$ is an $\omega$-compatible $S^1$-invariant almost complex structure. Let $\bigtriangledown H$ be the gradient vector field with respect to $g$, i.e., $dH = g(\bigtriangledown H,\cdot)$.
    Since $H$ is a Morse-Bott function by Corollary \ref{corollary : collary of equiv Darboux}, the unstable submanifold of $Z_{\min}$
    $$W^u(Z_{\min}) = \{p \in M ~|~ \lim_{s\rightarrow -\infty} \exp (s{\bigtriangledown H}) \cdot p \in Z_{\min} \}$$ is open and dense in $M$.
    Similarly, let $W^s(Z_{\max})$ be the stable submanifold of $Z_{\max}$. Since both $W^u(Z_{\min})$ and $W^s(Z_{\max})$ are open dense subsets, their intersection
    $W^u(Z_{\min}) \cap W^s(Z_{\max})$ is also open and dense in $M$, in particular it is non-empty.
    Now, pick a point $p \in W^u(Z_{\min}) \cap W^s(Z_{\max})$ and let
    \begin{displaymath}
        \begin{array}{ccc}
            M(p)& = & \displaystyle \bigcup_{t \in S^1, s \in \R} t \cdot \exp (s{\bigtriangledown H}) \cdot p. \\
        \end{array}
    \end{displaymath}
    Then it is straightforward that the closure $\overline{M(p)}$ is homeomorphic to a two-sphere whose north pole, respectively south pole, is in $Z_{\min}$, respectively  $Z_{\max}$.

    Now, let $\underline{X}$ be the fundamental vector field generated by the $S^1$-action.
    Since the $S^1$-action and the gradient flow are smooth, it is obvious that $M(p)$ is a smoothly embedded two-sphere punctured at the two poles $\{z_N, z_S\}$.
    In particular, a tangent space $T_qM(p)$ for every $q \in M(p)$ is generated by $\underline{X}_q$ and $J \underline{X}_q$ since $\underline{X}_q$ never vanishes on $M(p)$.
    Since we have chosen $J$ such that $\omega(\underline{X}_q, J \underline{X}_q) = g(\underline{X}_q, \underline{X}_q) > 0$, $M(p)$ is a smoothly embedded symplectic two sphere punctured at $\{z_N, z_S\}$. To make $\overline{M(p)}$ to be smooth near the poles, we perturb $g$ as follows.

    Let us consider an equivariant Darboux coordinate system $(\mathcal{U}_{\min}; z_1,z_2,\cdots,z_n)$ centered at $Z_{\min}$ such that
    $\omega = \frac{1}{2i} \sum dz_j \wedge d\bar{z_j}$ and the $S^1$-action is coordinate-wise linear.
    Let us take a standard almost complex structure $J_{\min}$ on $\mathcal{U}_{\min}$ with respect to the $(z_1,z_2,\cdots,z_n)$ coordinate system.
    Also, let $g$ be an $S^1$-invariant Riemannian metric on $\mathcal{U}_{\min}$ defined by $g_{\min}(X,Y) = \omega(X,J_{\min}Y)$ on $\mathcal{U}_{\min}$.
    Then it is easy to check that for any point $p \in \mathcal{U}_{\min}$, $\overline{M(p)}$ is a complex one-dimensional subspace so that $\overline{M(p)}$ is smooth in $\mathcal{U}_{\min}$.
    Similarly, we can choose a standard $\omega$-compatible metric $g_{\max}$ near $Z_{\max}$. Then we may perturb $g$ near $Z_{\min}$ and $Z_{\max}$ by using an $S^1$-invariant bump functions which is $1$ near the extremum, and $0$ outside of $\mathcal{U}_{\min}$ and $\mathcal{U}_{\max}$. Then the perturbed metric $\widetilde{g}$ coincides with $g_{\min}$ near $Z_{\min}$ and $g_{\max}$ near $Z_{\max}$ respectively. Then $\overline{M(p)}$ with respect to $\widetilde{g}$ is smooth and symplectic.

    Now, suppose that $S = \overline{M(p)}$ be a smoothly embedded $S^1$-invariant symplectic two sphere which attains the maximum and the minimum as we constructed above.
    The restriction $TM|_S$ of the tangent bundle $TM$ to $S$ is an $S^1$-equivariant complex vector bundle over $S$. To complete the proof of Proposition \ref{proposition : embedding sphere Hamiltonian case}, we need the following lemma.

    \begin{lemma}\label{lemma : equivariant vector bundle over S^2}\cite{Au}
        Let $\pi : E \rightarrow S$ be an $S^1$-equivariant complex vector bundle over a two sphere $S$. Let $z_N$ ($z_S$, respectively) be the north (south, respectively) pole fixed by the action. Then we have $$\langle c_1(E), [S] \rangle = \frac{\lambda_{z_S} - \lambda_{z_N}}{k}$$ where $k$, respectively $-k$, is the rotation number near $z_N$, respectively near $z_S$,
        $\lambda_{z_N}$ ($\lambda_{z_S}$, resp.) is the sum of all weights of the $S^1$-representation over the fiber $\pi^{-1}(z_N)$ ($\pi^{-1}(z_S)$, resp.).
    \end{lemma}

   Since the sum of all weights $\lambda_{\max}$, respectively $\lambda_{\min}$, at $Z_{\max}$, respectively $Z_{\min}$, is negative, respectively positive, we have $$\langle c_1(TM), [S] \rangle = \frac{\lambda_{\min} - \lambda_{\max}}{k} > 0.$$

\end{proof}

\bigskip
\bibliographystyle{amsalpha}

\end{document}